Volodymyr Denysiuk

# The method of summation of divergent trigonometric series

## Abstract

The generalized summation of divergent trigonometric series, namely by method of $\sigma_k(r,\alpha)$- factors is considered in this paper. It is proved that such summation of Fourier series of periodical function $f(t)$ results in the convolution of this function with kernels $De(r,\alpha,t)$; if the parameter $r$ is integer, these kernels are polynomial normalized basic $B$-splines of order $r$-1 ($r = 1, 2, ...$). Also it is proved that the method of summation with $\sigma_k(r,\alpha)$- factors is $F$-effective.

## Introduction

Consider a divergent trigonometric series

$$\frac{a_0}{2} + \sum_{k=1}^{\infty}\left[a_k \cos kt + b_k \sin k\alpha\right], \qquad (1)$$

with coefficients $|a_k|, |b_k| \leq C$, ($0 < C < \infty$, $k = 1, 2, ...$), and $\lim_{k\to\infty}(a_k, b_k) \neq 0$.

To provide uniform convergence of series of type (1), generalized methods of summation of divergent series are often applied. Among these methods, we consider some generalized summation methods of the Poisson-Abel type. As is known [1], the basic idea of these summation methods of a series (1) is that all terms of a series are simultaneously multiplied by factors of a certain type; in other words, the series (1) corresponds to the series of the form:

$$\frac{a_0}{2} + \sum_{k=1}^{\infty}\left[a_k \cos kt + b_k \sin k\alpha\right] \to \frac{a_0}{2} + \sum_{k=1}^{\infty}\mu_k(\alpha_1, \alpha_2, \ldots, \alpha_n)\left[a_k \cos kt + b_k \sin k\alpha\right], \quad (2)$$

where $\mu_k(\alpha_1, \alpha_2, \ldots, \alpha_n)$, ($k = 1, 2, ...$) are some factors dependent on $n$ ($n = 1, 2, \ldots$) parameters $\alpha_1, \alpha_2, \ldots, \alpha_n$; such factors are often called convergence factors. Accordingly, the function $f(t)$ corresponds to function $f(\alpha_1, \alpha_2, \ldots, \alpha_n, t)$, where

$$f(\alpha_1, \alpha_2, \ldots, \alpha_n, t) \square \frac{a_0}{2} + \sum_{k=1}^{\infty}\mu_k(\alpha_1, \alpha_2, \ldots, \alpha_n)\left[a_k \cos kt + b_k \sin k\alpha\right]. \qquad (3)$$

If the series (3) converges uniformly, then it is a Fourier series of the function $f(\alpha_1, \alpha_2, \ldots, \alpha_n, t)$, and the values

$$a_0, \ \mu_k(\alpha_1, \alpha_2, \ldots, \alpha_n) a_k, \ \mu_k(\alpha_1, \alpha_2, \ldots, \alpha_n) b_k, \ (k=1,2,\ldots)$$

are Fourier coefficients of this function. This fact deserves attention, since the coefficients $a_0$, $a_k$, $b_k$, ($k=1,2,\ldots$) of the series (1) have no relation to the Fourier coefficients.

According to the type of the factors $\mu_k(\alpha_1, \alpha_2, \ldots, \alpha_n)$ specific methods of summation are distinguished. We will restrict our consideration to the Poisson-Abel method and the $(R,k)$-method [1] which we will call method of $\sigma_k(r,\alpha)$-factors.

The Poisson-Abel method suggests that $\mu_k(\alpha_1, \alpha_2, \ldots, \alpha_n) = r^k$, (i.e. $n=1$, $\alpha_1 = r$). According to (2), the function $f(t)$ (1) corresponds to series

$$f(r,t) \square \frac{a_0}{2} + \sum_{k=1}^{\infty} r^k \left[ a_k \cos kt + b_k \sin k\alpha \right], \ (0 < r < 1). \tag{4}$$

It is obvious that the series (4) is uniformly convergent, since it is majorized by a geometric progression $C \sum_{k=1}^{\infty} r^k$. This derives that the coefficients $a_0, r^k a_k, r^k b_k$, ($k=1,2,\ldots$), are Fourier coefficients of the function $f(r,t)$.

It would be convenient to represent the series in (4) as a convolution. In order to get such a representation, we will consider the following series:

$$\frac{1}{2} + \sum_{k=1}^{\infty} r^k \cos k(u-t), \quad (0 < r < 1) \tag{5}$$

which is uniformly convergent.

By means of term-by-term multiplication of the series (5) by the absolutely integrable function $f(t)$, with simple transformations and term-by-term integration we get

$$f(r,t) \square \frac{a_0}{2} + \sum_{k=1}^{\infty} r^k \left[ a_k \cos kt + b_k \sin k\alpha \right] = \frac{1}{2\pi} \int_{-\pi}^{\pi} f(u+t) \left\{ \frac{1}{2} + \sum_{k=1}^{\infty} r^k \cos ku \right\} du. \tag{6}$$

To obtain (6) we used the fact that the uniformly convergent series

$$\frac{1}{2} + \sum_{k=1}^{\infty} r^k \cos k(u-t), \quad (0 < r < 1),$$

all terms of which have been multiplied by an absolutely integrable function, may be integrated termwise.
Considering that

$$\frac{1}{2} + \sum_{k=1}^{\infty} r^k \cos ku = \frac{1-r^2}{1 - 2r \cos u + r^2},$$

we eventually get

$$f(r,t) = \frac{a_0}{2} + \sum_{k=1}^{\infty} r^k \left[ a_k \cos kt + b_k \sin k\alpha \right] = \frac{1}{2\pi} \int_{-\pi}^{\pi} f(u+t) P(r,u) \, du. \tag{7}$$

This integral is referred to as the Poisson integral, and the kernel $P(r,u)$ of the integral transformation

$$P(r,u) = \frac{1}{2} \frac{1-r^2}{1-2r\cos u + r^2}$$

as the Poisson kernel. The Poisson integral has been theoretically justified by Schwarz.

So, (7) can be presented as

$$f(r,t) = \frac{a_0}{2} + \sum_{k=1}^{\infty} r^k \left[ a_k \cos kt + b_k \sin k\alpha \right] = \frac{1}{2\pi} \int_{-\pi}^{\pi} f(u+t) P(r,u) du.$$

It is easy to verify that with $r \to 1-0$, the Poisson kernel generates a $\delta$ - like series [8]. From this follows a convergence theorem.

**Theorem.** At the point where a function $f(t)$ is continuous or, at least, it has a discontinuity of the first kind, the Fourier series (1) can be summed by the Poisson-Abel method, with $.5[f(t-0) + f(t+0)]$ being the generalized sum of the series. It is well understood that at the continuity point the sum of the series is equal to $f(t)$.

Thus, we can conclude that the Poisson-Abel summation method is $F$ - effective [1].

A disadvantage of the Poisson-Abel method is that the parameter $r$ has no interpretation. Therefore only the limit function $f(1-0,t) = \lim_{r \to 1-0} f(r,t)$ is taken into consideration when this method is applied. However, in many cases limit transition is not practical. Thus attention is drawn to methods in which parameters defining the method may be interpreted in one way or another. Applying such methods, we can consider not only limit functions obtained as a result of limit transition by some parameters, but also functions derived with these parameters being fixed. One of such methods is the summation method with $\sigma_k(r,\alpha)$ factors, which we consider here.

**The aim of the work**

The aim of the work consists in studying the method of Fourier trigonometric series summation with $\sigma_k(r,\alpha)$ - factors.

**Main body**

In the summation method with $\sigma_k(r,\alpha)$ - factors we assume

$$\mu_k(\alpha_1, \alpha_2, \ldots, \alpha_n) = \sigma_k(r,\alpha),$$

where

$$\sigma_k(r,\alpha) = \left( \frac{\sin k\frac{\alpha}{2}}{k\frac{\alpha}{2}} \right)^r, \qquad (0 < \alpha < \pi;\ r = 1,2,\ldots). \qquad (8)$$

It is clear that in this case $n=2$, $\alpha_1 = \alpha$; $\alpha_2 = r$.

Further for convenience we will use the following notation:

$$\operatorname{sinc} x = \frac{\sin x}{x}.$$

In this notation the expression for $\sigma_k(r,\alpha)$ becomes simpler:

$$\sigma_k(r,\alpha) = \left(\operatorname{sinc} k\frac{\alpha}{2}\right)^r.$$

To begin with, we will note that the factors $\sigma_k(2,2h)$ were derived by Riemann when performing formal double integration of trigonometric series with further determination of Schwarz's generalized second derivative [2]. The factors were used to sum trigonometric series with further limit transition $h \to 0$.

The factors $\sigma_k(1,\alpha)$ may also occur as a result of application of the method of phantom functions [3].

Factors of the type $\sigma_k(1,\frac{2\pi}{n+1})$ ($k=1,2,...,n$) were derived by K. Lanczos when smoothing the partial sums $S_n$ of a Fourier series [4]; this factor was used in $\lambda$-summation methods. Lastly, powers of such factors may also occur when determining Fourier approximate coefficients by the Filon method [5].

The summation method with $\sigma_k(r,\alpha)$ - factors suggests that the series (1) has its corresponds to the series

$$\frac{a_0}{2} + \sum_{k=1}^{\infty} \sigma_k(r,\alpha)\left[a_k \cos kt + b_k \sin k\alpha\right], \quad (r \geq 1) \qquad (9)$$

where $\sigma_k(r,\alpha)$ are defined by (8).

It is well understood that with $r>1$ the series (9) is uniformly convergent. Indeed, as the coefficients $a_k, b_k$ of the series tend to 0 with $k \to \infty$, they are collectively bounded

$$|a_k|, |b_k| < K, \quad (K = const),$$

and the series (6) is majorized by the convergent series $2K\sum_{k=1}^{\infty}\frac{1}{k^r}$. Therefore, according to the Weierstrass criterion, the series (9) is uniformly convergent.

As in the case before, it is easy to derive the expression

$$f(r,\alpha,t) = \frac{1}{2\pi} \int_{-\pi}^{\pi} f(u+t)\left\{\frac{1}{2} + \sum_{k=1}^{\infty} \sigma_k(r,\alpha)\cos ku\right\} du. \qquad (10)$$

Introduce the notation

$$De(r,\alpha,t) = \frac{1}{2} + \sum_{k=1}^{\infty} \sigma_k(r,h)\cos k(t) \qquad (11)$$

The kernel $De(r,\alpha,t)$ will be referred to as De-kernel of the $r$-th order. Consider these kernels in more detail, noting the following.

Let $r=1$.

First of all, it is easy to see that the series which defines the kernel $De(1,\alpha,t)$ is convergent (but not uniformly convergent) according to the Dirichlet criterion. Determining the sum of this series, we get:

$$De(1,\alpha,t) = \frac{1}{\pi}\left(\frac{1}{2} + \sum_{k=1}^{\infty}\left(\text{sinc}\,k\frac{\alpha}{2}\right)\cos kt\right) =$$

$$= \frac{1}{\pi}\left[\frac{1}{2} + \frac{1}{\alpha}\sum_{k=1}^{\infty}\frac{\sin k\left(t+\frac{\alpha}{2}\right)}{k} - \frac{1}{\alpha}\sum_{m=1}^{n}\frac{\sin k\left(t-\frac{\alpha}{2}\right)}{k}\right]$$

With simple transformations we get

$$De(1,\alpha,t) = \begin{cases} \dfrac{1}{\alpha}, & |t| < \dfrac{\alpha}{2} \\ 0, & \dfrac{\alpha}{2} < |t| \leq \pi \end{cases}.$$

Thus, the kernel $De(1,\alpha,t)$ is a finite discontinuous function.

Now consider the kernel $De(2,\alpha,t)$. One can immediately verify that this kernel can be represented as

$$De(2,\alpha,t) = \frac{1}{\pi}\left(\frac{1}{2} + \sum_{k=1}^{\infty}\left(\text{sinc}\,k\frac{\alpha}{2}\right)^2\cos kt\right) = \begin{cases} \dfrac{\alpha-|t|}{\alpha^2}, & |t| \leq \alpha \\ 0, & \alpha < |t| \leq \pi \end{cases}.$$

It is clear that the trigonometric series defining this kernel is uniformly convergent.

As before, the properties of the kernel imply that

$$\int_{-\pi}^{\pi} De(2,\alpha,t)\,dt = 1.$$

Now consider an explicit representation of the kernel $De(r,\alpha,t)$ ($r = 3, 4, ...$). One can immediately verify that the kernel $De(r,\alpha,t)$, up to the constant factor $\alpha$, coincides with periodically continued polynomial normalized $B$-splines of the order $r-1$ ($r = 1, 2, ...$) with $\alpha = h$, which are constructed on uniformly spaced grids with a step $h$, see [6], [7]. Denoting such splines as $B_{r-1}(\alpha,t)$ we get

$$\frac{1}{\alpha}B_{r-1}(\alpha,t) = \frac{1}{\pi}\left(\frac{1}{2} + \sum_{k=1}^{\infty}\left(\text{sinc}\,\frac{k\alpha}{2}\right)^r\cos kt\right), \qquad (12)$$

where the spline $B_{r-1}(\alpha,t)$ is constructed on a uniformly spaced grid with the step $\alpha$, the symmetry centers of both functions being concordant. We will note that for $r = 1,...,5$ explicit representation of $B$-splines can be used [5].

Due to the properties of $B$-splines (12) can be presented as

$$B_r(\alpha,t) = \frac{1}{\alpha}\int\limits_{-\alpha/2}^{\alpha/2} B_{r-1}(\alpha,t-u)B_0(\alpha,u)du = \quad (13)$$

$$= \frac{1}{\pi}\int\limits_{-\alpha/2}^{\alpha/2}\left(\frac{1}{2} + \sum_{k=1}^{\infty}\left(\frac{\sin\frac{k\alpha}{2}}{\frac{k\alpha}{2}}\right)^r \cos k(t-u)\right)B_0(\alpha,u)du = \frac{\alpha}{\pi}\left(\frac{1}{2} + \sum_{k=1}^{\infty}\left(\frac{\sin\frac{k\alpha}{2}}{\frac{k\alpha}{2}}\right)^{r+1} \cos kt\right).$$

In our view, expression (13) deserves special attention as it establishes relationship between convolution-type integrals for $B$-splines and their representation in terms of a Fourier series.

Thus, considering (7), the function $f(r,\alpha,t)$ can be presented as

$$f(r,\alpha,t) = \frac{1}{\alpha}\int\limits_{-\pi}^{\pi} f(t+u)B_{r-1}(\alpha,u)\,du, \qquad r = 1, 2, \ldots.$$

It is easy to verify that the kernels $De(r,\alpha,t)$, for any fixed $r = 1, 2, \ldots$, form a $\delta$-like series [8] as $\alpha \to 0$. From this follows a convergence theorem.

**Theorem.** At the point where a function $f(t)$ is continuous or, at least, it has a discontinuity of the first kind, the Fourier series (1) can be summed by the method of of $\sigma_k(r,\alpha)$- factors, with $.5[f(t-0) + f(t+0)]$ being the generalized sum of the series. It is well understood that at the continuity point the sum of the series is equal to $f(t)$.

$$\lim_{\alpha\to 0} f(r,\alpha,t) = \lim_{\alpha\to 0}\frac{1}{\alpha}\int\limits_{-\pi}^{\pi} f(t+u)B_{r-1}(\alpha,u)\,du = \frac{1}{2}[f(t-0) + f(t+0)].$$

In particular, in the continuity point this limit equals to $f(t)$.

So, same as above, we may conclude that the summation method with $\sigma_k(r,\alpha)$-factors is $F$-effective [24].

It is well understood that the parameters $r$ and $\alpha$, which are included in the expression for $\sigma_k(r,\alpha)$-factors, allow for an easy interpretation. In particular, the parameter $r$ defines the differential features of the kernel $De(r,\alpha,t)$ whereas the parameter $\alpha$, for functions that have points of discontinuity of the first kind (jump-type discontinuity), defines the neighborhoods of these points in which these discontinuities are smoothed.

**Conclusions**

In this paper we consider and study a generalized summation method of trigonometric series with $\sigma_k(r,\alpha)$-factors. Here we show that this method leads to the kernels $De(r,\alpha,t)$ of integral transformation. It is also shown that, for integer values of parameter $r$, these kernels are polynomial normalized basis $B$-splines of the corresponding orders.

It is established that the sum of the Fourier series of an absolutely integrable function $f(t)$ with introduced convergence factors $\gamma_k(\alpha_1,\alpha_2,\ldots,\alpha_n)$ provided that the series is uniformly convergent is equal to the convolution of this function with the kernel $De(t,\alpha_1,\alpha_2,\ldots,\alpha_n)$, which is determined by the equation

$$De(t,\alpha_1,\alpha_2,\ldots,\alpha_n) = \frac{1}{\pi}\left(\frac{1}{2} + \sum_{k=1}^{\infty}\gamma(\alpha_1,\alpha_2,\ldots,\alpha_n)\cos kt\right).$$

In particular, it is shown here that the trigonometric Fourier series of an absolutely integrable function $f(t)$, which are summed by the summation method with $\sigma_k(r,\alpha)$-factors, converge to the convolution of this function with polynomial $B$- splines. The establishment of this relationship offers great opportunities to research the classes of polynomial and trigonometric approximations as a whole.

The factors $\sigma_k(r,\alpha)$ are the Fourier coefficients of normalized basis $B$-splines of the order $r-1$, ($r=1,2,\ldots$), whereas the kernels $De(r,\alpha,t)$ may be considered a trigonometric representation of such splines. As the $B$-splines form a polynomial basis in the spaces $C^r$, we can assert that the kernels $De(r,\alpha,t)$ also form a basis in these spaces. Therefore it is possible to construct trigonometric analogues of simple polynomial splines via linear combinations of these kernels.

Expression (13), which establishes the relationship between convolution-type integrals for $B$-splines and their representation in terms of a Fourier series, also attracts attention.

There is an opportunity to use other types of finite functions with a fixed smoothness, in particular spline-type functions, as $De(r,\alpha,t)$ kernels. Moreover, it becomes possible to construct kernels of the $De(r,\alpha,t)$ type with given properties. The summation factors of the $\sigma_k(r,\alpha)$ type for such kernels can be easily obtained as Fourier coefficients in expansions of these kernels according to (12).

Lastly, one more reason why the trigonometric functions summation method with $\sigma_k(r,\alpha)$ factors attracts attention is that it represents a generalization of the results described in [1,2,9]. For example, A. Zygmund and G. Hardy represent the $De(2,\alpha,t)$ kernel in terms of $\sigma_k(2,\alpha)$ factors, and N. Bari – the $De(3,\alpha,t)$ kernel in terms of $\sigma_k(3,\alpha)$ factors. Nevertheless, these authors did not go any further towards the generalization of these results.

**List of references**

1. Divergent Series by G. H Hardy, Oxford, 1949.
2. Bari N.K. Trigonometric Series. – M.: Gosizdat fiz.-mat. lit., 1961. – 936 p. (Russian).
3. Denysiuk V.P. Method of improving of convergence for Fourier series and interpolating polynomial in orthogonal function // Naukovi visti NTTU "KPI". – 2013. - № 4. – P.31 – 37. (Ukrainian).
4. Lanchos C. Applied Analysis. Prentice Hall, Inc., 1956.

5. Denysiuk V.P. Trigonometric series and splines. - Kyiv, NAU, 2017 . – 212 p. (Ukrainian).

6. Zav'yalov Yu.S., Kvasov B. I., Miroshnichenko V. L. Methods of spline-functions. – M.: Nauka, 1980. – 352 p. (Russian).

7. de Boor C. A Practical Guide to Splines. - M.: Radio i sviaz, 1985. – 304 p. (Russian).

8. Dziadyk V. K. Introduction to theory of uniform functions approximation by polynomials . – M.: Nauka, 1977. – 512 p. (Russian).

9. Zygmund A. Trigonometric series.- Volume 1,2, Cambridge at The University Press, 1959.